\documentclass[10pt]{amsart}

\usepackage{latexsym,enumerate}
\usepackage{amsmath,amsthm,amsopn,amstext,amscd,amsfonts,amssymb}
\usepackage[ansinew]{inputenc}
\usepackage{verbatim}
\usepackage{wasysym}
\usepackage[all]{xy}
\newcommand{\CC}{\mathbb{C}}
\newcommand{\NN}{\mathbb{N}}

\newcommand{\RR}{\mathbb{R}}

\newcommand{\lin}{\mathcal{L}}

\newtheorem{theorem}{Theorem}[section]

\newtheorem{proposition}[theorem]{Proposition}
\newtheorem{corollary}[theorem]{Corollary}

\DeclareMathOperator{\diam}{diam}

\DeclareMathOperator{\Deg}{Deg}

\renewcommand{\a}{\alpha}

\renewcommand{\d}{\delta}
\newcommand{\D}{\Delta}
\newcommand{\g}{\gamma}
\newcommand{\G}{\Gamma}

\begin{document}
%\DeclareGraphicsExtensions{.jpg,.pdf,.mps,.png}
\title{Mathematical results on harmonic polynomials}

\author[W. Carballosa]{Walter Carballosa}
\address{Department of Mathematics and Statistics, Florida International University, 11200 SW 8th Street
Miami, FL 33199, USA.}
\email{waltercarb@gmail.com}

\author[J. E. N\'apoles]{Juan E. N\'apoles}
\address{Departamento de Matem\'aticas, Universidad Nacional de Nordeste, Avenida de la Libertad 5450, 3400 Corrientes, Argentina}
\email{jnapoles@frre.utn.edu.ar}
%\thanks{$^{(1)}$ Supported in part by a grant from CONACYT (FOMIX-CONACyT-UAGro 249818), M\'exico.}

\author[J. M. Rodr{\'\i}guez]{Jos\'e M. Rodr{\'\i}guez$^{(1)}$}
\address{Departamento de Matem\'aticas, Universidad Carlos III de Madrid, Avenida de la Universidad 30, 28911 Legan\'es, Madrid, Spain}
\email{jomaro@math.uc3m.es}
\thanks{$^{(1)}$ Supported in part by two grants from Ministerio de Econom{\'\i}a y Competititvidad, Agencia Estatal de
Investigación (AEI) and Fondo Europeo de Desarrollo Regional (FEDER) (MTM2016-78227-C2-1-P and MTM2017-90584-REDT), Spain.}

\author[O. Rosario]{Omar Rosario}
\address{Facultad de Matem\'aticas, Universidad Aut\'onoma de Guerrero, Carlos E. Adame No.54 Col. Garita, 39650 Acalpulco Gro., Mexico}
\email{omarrosarioc@gmail.com}

\author[J. M. Sigarreta]{Jos\'e M. Sigarreta$^{(1)}$}
\address{Facultad de Matem\'aticas, Universidad Aut\'onoma de Guerrero, Carlos E. Adame No.54 Col. Garita, 39650 Acalpulco Gro., Mexico}
\email{jsmathguerrero@gmail.com}

\date{\today}

\maketitle{}

%\centerline{{\bf HARMONIC polynomial 21.TEX}}

\begin{abstract}
Some years ago, the harmonic polynomial was introduced in order to understand better the harmonic topological index; for instance,
it allows to obtain bounds of the harmonic index of the main products of graphs.
Here, we obtain several properties of this polynomial, and we prove that several properties of graphs can be deduced from their harmonic polynomials.
Also, we show that two graphs with the same harmonic polynomial have to be similar.
\end{abstract}



\section{Introduction}

A topological descriptor is a single number that represents a chemical structure in graph-theoretical terms via the
molecular graph, they play a significant role in mathematical chemistry especially in the QSPR/QSAR investigations.
A topological descriptor is called a topological index if it correlates with a molecular property.
Topological indices are used to understand physicochemical properties of chemical compounds,
since they capture some properties of a molecule in a single number.
Hundreds of topological indices have been introduced and studied, starting with the
seminal work by Wiener \cite{Wi}.

Within all topological indices ones of the
most investigated are the descriptors based on the valences of atoms in molecules (in
graph-theoretical notions degrees of vertices of graph).
Among them, several indices are recognized to be useful tools in
chemical researches. Probably, the best know such descriptor is the Randi\'c connectivity
index ($R$) \cite{R}. There are more than thousand papers and a couple of books dealing with
this molecular descriptor (see, e.g., \cite{GF}, \cite{LG}, \cite{LS}, \cite{RS}, \cite{RS0} and the references therein).
During many years, scientists were trying to improve the predictive power of the
Randi\'c index. This led to the introduction of a large number of new topological
descriptors resembling the original Randi\'c index.
Two of the main successors of the Randi\'c index are the first and second Zagreb indices,
denoted by $M_1$ and $M_2$, respectively, defined as
$$
M_1(G) = \sum_{uv\in E(G)} (d_u + d_v) = \sum_{u\in V(G)} d_u^2,
\qquad
M_2(G) = \sum_{uv\in E(G)} d_u d_v ,
\qquad
$$
where $uv$ denotes the edge of the graph $G$ connecting the vertices $u$ and $v$, and
$d_u$ is the degree of the vertex $u$.
These indices have attracted growing interest, see e.g. \cite{BF}, \cite{D2}, \cite{FGE}, \cite{L}
(in particular, they are included in a number of programs used for the routine computation of topological indices).

Another remarkable topological descriptor is the \emph{harmonic} index, defined in \cite{Faj} as
$$
H(G) = \sum_{uv\in E(G)}\frac{2}{d_u + d_v}\,.
$$
This index has attracted a great interest in the lasts years
(see, e.g., \cite{DBAV}, \cite{FMS}, \cite{RS harm}, \cite{WTD}, %\cite{Z}
and \cite{ZX}).

With motivation from the first Zagreb and harmonic indices, \emph{general sum-connectivity index} $\chi_{\a}$
was defined by Zhou and Trinajsti\'c in \cite{ZT2} as
$$
\chi_{\a}(G) = \sum_{uv\in E(G)} (d_u + d_v)^{\a},
$$
with $\a \in \RR$.
Note that $\chi_{_{1}}$ is the first Zagreb index $M_1$, $2\chi_{_{-1}}$ is the harmonic index $H$,
$\chi_{_{-1/2}}$ is the sum-connectivity index, etc.
Some mathematical properties of the general sum-connectivity index were given in
\cite{DZT3}, \cite{RS harm}, \cite{ZX}, \cite{ZT2} and \cite{ZL}.

\smallskip

In \cite{Iranmanesh Saheli} the \emph{harmonic polynomial} of a graph $G$ is defined as
$$
H(G,x) := \sum_{uv\in E(G)} x^{d_u + d_v - 1},
$$
and the harmonic polynomials of some graphs are computed.
The harmonic polynomials of the line of some graphs are computed in \cite{Nazir Sardar Zafar Zahid}.
In \cite{HMRS} this polynomial was used in order to obtain bounds of the harmonic index of the main products of graphs.

This polynomial gets its name from the fact that $2\int^1_0 H(G, x) \,dx = H(G)$.

If $G_1$ and $G_2$ are disjoint graphs, then
$$
H(G_1\cup G_2, x) = H(G_1, x) + H(G_2, x).
$$
Hence, considering connected graphs is not a restrictive condition.

\smallskip

The characterization of any graph by a polynomial is one of the open important problems in graph theory. In recent years there have been many works on graph polynomials (see, e.g., \cite{book} and the references therein). %(see, e.g., \cite{Alk,carb1,carb2,Graph Polynomials,tute}).
The research in this area has been largely driven by the advantages offered by the use of computers: it is simpler to represent a graph by a polynomial (a vector with dimension $O(n)$) than by the adjacency matrix (an $n\times n$ matrix). Some parameters of a graph allow to define polynomials related to a graph. Although several polynomials are interesting since they compress information about the graph’s structure, unfortunately, the well-known polynomials do not solve the problem
of the characterization of any graph, since there are often non-isomorphic graphs with the same polynomial.

\smallskip

Throughout this paper, $G=(V,E)=(V (G),E (G))$ denotes a (non-oriented) finite simple (without multiple edges and loops)
graph without isolated vertices (every vertex has at least a neighbor).
The main aim of this paper is to obtain several properties of the harmonic polynomial.
We prove that several properties of graphs can be obtained from their harmonic polynomials:
Corollary \ref{c:p2}
characterizes regular and biregular graphs in terms of the zeros of their harmonic polynomials;
Theorem \ref{t:p3} gives information about the connectedness, the diameter and the girth (the minimum length of the cycles) of a graph
in terms of the degree of its harmonic polynomial;
Proposition \ref{p:5} shows that the cardinality of the set of pendant paths in a graph is precisely the coefficient of $x^2$ in its harmonic polynomial.
Besides, Theorems \ref{t:sucdeg}, \ref{t:sucdeg2} and \ref{t:sucdeg3} relate the number of non-zero coefficients of the harmonic polynomial with the degree sequence of the polynomial.
Theorem \ref{t:mu} shows that two graphs with the same harmonic polynomial have to be similar.

\section{Main results}

The following result appears in \cite[Proposition 1]{Iranmanesh Saheli}.

\begin{proposition} \label{Poli1}
If $G$ is a $k$-regular graph with $m$ edges, then
$H(G,x)=mx^{2k-1}.$
\end{proposition}

Propositions 2, 4, 5, 7
in \cite{Iranmanesh Saheli}
have the following consequences on the graphs:
$K_{n}$ (the complete graph with $n$ vertices),
$C_{n}$ (the cycle with $n\ge 3$ vertices),
$Q_{n}$ (the $n$-dimensional hypercube),
$K_{n_1,n_2}$ (the complete bipartite graph with $n_1+n_2$ vertices),
$P_{n}$ (the path graph with $n$ vertices) and
$W_{n}$ (the wheel graph with $n\geq4$ vertices).

\begin{proposition} \label{Poli2}
We have
$$
\begin{aligned}
H(K_{n},x)
 =\frac12\,n(n-1)x^{2n-3},
\qquad
& H(C_{n},x)
 =nx^{3},
\\
H(Q_{n},x)
= n2^{n-1}x^{2n-1},
\qquad
& H(K_{n_1,n_2},x)
=n_1n_2x^{n_1+n_2-1},
\\
H(P_{n},x) =2x^{2}+(n-3)x^{3},
\qquad
& H(W_{n},x)
=(n-1)(x^{n+1}+x^{5}).
\end{aligned}
$$
\end{proposition}

Given a graph $G$, the \emph{line graph} $\lin(G)$ of $G$ is a graph which has a vertex $w_{e}\in V(\lin(G))$
for each edge $e \in E(G)$, and an edge joining $w_{e_i}$ and $w_{e_j}$ when $e_i$ and $e_j$ share a vertex
(i.e., $\lin(G)$ is the intersection graph of $E(G)$).
It is easy to check that if $uv \in E(G)$, then the degree of $w_{uv} \in V(\lin(G))$ is $d_u + d_v - 2$.

Line graphs were initially introduced in the papers \cite{W} and \cite{Kr},
although the terminology of line graph was used in \cite{HN} for the first time.
They are an active topic of research at this moment.

In the same paper, where Zagreb indices were introduced, the \emph{forgotten topological index} (or \emph{F-index}) is defined as
$$
F(G) = \sum_{uv\in E(G)} (d_u^2 + d_v^2) = \sum_{u \in V(G)} d_u^3.
$$
Both the forgotten
topological index and the first Zagreb index were employed in the formulas for total $\pi$-electron energy in \cite{GT},
as a measure of branching extent of the carbon-atom skeleton of the underlying molecule.
However, this index never got attention
except recently, when Furtula and Gutman in \cite{3} established some basic properties
of the F-index and showed that its predictive ability is almost similar to
that of first Zagreb index and for the entropy and acetic factor, both of them yield correlation
coefficients greater than $0.95$.

Our first result shows that we can obtain information about the graph from the values of the harmonic polynomial (and its derivatives) at the point $1$.

\begin{proposition} \label{p:1}
If $G$ is a graph with $n$ vertices, $m$ edges, maximum degree $\D$ and minimum degree $\d$, then:
\begin{itemize}
\item $H(G,1) = m$,
\item $H'(G,1) + H(G,1)= M_1(G)$,
\item $H''(G,1)-2H(G,1) = F(G)+2M_2(G)-3M_1(G)$,
\item $H''(G,1) +2H(G,1)= M_1(\lin (G)) + M_1(G)$,
\item $2H(G,1)/\D \le n \le 2H(G,1)/\d$.
\end{itemize}
\end{proposition}

\begin{proof}
First of all, $H(G,1)=\sum_{uv\in E(G)} 1 =m$.
Also,
$$
\begin{aligned}
%H'(G,x) & = \sum_{uv\in E(G)} \!\!\!(d_u + d_v - 1) \, x^{d_u + d_v - 2},
%\\
H'(G,1) & = \sum_{uv\in E(G)} \!\!\!(d_u + d_v) -\!\!\! \sum_{uv\in E(G)} \!\!\! 1
= M_1(G) - H(G,1),
\end{aligned}
$$
and
$$
\begin{aligned}
%H''(G,x) & = \sum_{uv\in E(G)} \!\!\!(d_u + d_v - 1)(d_u + d_v - 2) \, x^{d_u + d_v - 3},
%\\
H''(G,1) & = \sum_{uv\in E(G)} \!\!\!(d_u^2 + d_v^2) +2 \!\!\!\sum_{uv\in E(G)} \!\!\!d_u d_v - 3\!\!\!\sum_{uv\in E(G)} \!\!\!(d_u + d_v) + \!\!\!\sum_{uv\in E(G)} \!\!\!2
\\
& = F(G)+2M_2(G)-3M_1(G)+2H(G,1),
\\
H''(G,1) & = \sum_{uv\in E(G)} \!\!\!(d_u + d_v - 2)(d_u + d_v - 2) + \sum_{uv\in E(G)} \!\!\!(d_u + d_v - 2)
\\
& = \sum_{uv\in E(G)} \!\!\!(d_u + d_v - 2)^2 + \sum_{uv\in E(G)} \!\!\!(d_u + d_v) - \!\!\!\sum_{uv\in E(G)} \!\!\!2
\\
& = M_1(\lin (G)) + M_1(G) - 2H(G,1).
\end{aligned}
$$
\indent
The inequalities $\d n \le 2m \le \D n$ and the first item imply the fifth one.
\end{proof}

Proposition \ref{Poli1} shows that any two $k$-regular graphs with the same cardinality of edges, have the same harmonic polynomial.
It is natural to ask the following question: How many graphs can be characterized by their harmonic polynomials?
This is a very difficult question, but there are partial answers:
Proposition \ref{p:1} gives that graphs with different cardinality of edges have different harmonic polynomials.
This fact has the following interesting consequence.

\begin{corollary} \label{c:1}
If $\G$ is a proper subgraph of the graph $G$, then
$H(\G,x) \neq H(G,x)$.
\end{corollary}

Also, Theorem \ref{t:mu} will show that two graphs with the same harmonic polynomial have to be similar, in some sense.

\smallskip

For each positive integer  $k$, let us define the polynomial
$$
Q_k(x)
:= (x-1)(x-2)\cdots(x-k)
= x^k + \sum_{j=0}^{k-1} a_{k,j}x^{j}.
$$
Note that Vieta's formulas allow to compute these coefficients $a_{k,j}$ in a very simple way:
$$
a_{k,k-j} = (-1)^j \!\!\!\!\!\!\!\!\!\!\!\! \sum_{1 \le i_1 < i_2 < \cdots < i_{j} \le k} \!\!\!\!\!\!\!\!\!\!\!\! i_1 i_2 \cdots i_{j}.
$$
%$$
%\begin{aligned}
%1+2+3+ \dots +k & = -a_{k,k-1} = \frac12\, k(k+1),
%\\
%(1\cdot 2+1\cdot 3 +  \dots +1\cdot k) + (2\cdot 3+2\cdot 4 +  \dots +1\cdot k) + \dots +(k-1)k & = a_{k,k-2},
%\\
%\dots
%\\
%1\cdot 2\cdot 3 \dots k & = (-1)^k a_{k,0}= k! .
%\end{aligned}
%$$
In particular, we have $a_{k,k-1} = - \frac12\, k(k+1)$ and $a_{k,0}= (-1)^k k!$.

\begin{proposition} \label{p:1k}
If $G$ is a graph and $k$ is a positive integer, then
$$
H^{(k)}(G,1)
= \chi_{k}(G) + \sum_{j=0}^{k-1} a_{k,j}\chi_{j}(G).
$$
\end{proposition}

\begin{proof}
We have
$$
\begin{aligned}
%H'(G,x) & = \sum_{uv\in E(G)} \!\!\!(d_u + d_v - 1) \, x^{d_u + d_v - 2},
%\\
%H''(G,x) & = \sum_{uv\in E(G)} \!\!\!(d_u + d_v - 1)(d_u + d_v - 2) \, x^{d_u + d_v - 3},
%\\
H^{(k)}(G,x) & = \sum_{uv\in E(G)} \!\!\!(d_u + d_v - 1)(d_u + d_v - 2)\cdots (d_u + d_v - k) \, x^{d_u + d_v - k-1}
\\
& = \sum_{uv\in E(G)} \!\!\!Q_k(d_u + d_v) \, x^{d_u + d_v - k-1},
\\
H^{(k)}(G,1)
& = \sum_{uv\in E(G)} \!\!\!Q_k(d_u + d_v)
= \sum_{uv\in E(G)} \!\!\!(d_u + d_v)^k + \sum_{j=0}^{k-1} \sum_{uv\in E(G)} \!\!\!a_{k,j}(d_u + d_v)^{j}
\\
& = \chi_{k}(G) + \sum_{j=0}^{k-1} a_{k,j}\chi_{j}(G).
\end{aligned}
$$
\end{proof}

As usual, we denote by $\Deg p(x)$ the degree of the polynomial $p(x)$,
and by $\Deg_{\min} p(x)$ the minimum degree of their monomials with non-zero coefficients.

\smallskip

Given a graph $G$, we have

$$
\begin{array}{c}
\Deg H(G,x) = \max\big\{ d_u+d_v-1\,|\; uv\in E(G) \big\}, \\ 
\Deg_{\min} H(G,x) = \min\big\{ d_u+d_v-1\,|\; uv\in E(G) \big\}.
\end{array} $$

Recall that a \emph{biregular} graph is a bipartite graph for which any vertex in one side of the given bipartition has degree $\D$
and any vertex in the other side of the bipartition has degree $\d$.
We say that a graph is $(\D,\d)$-\emph{biregular} if we want to write explicitly the maximum and minimum degrees.

\begin{proposition} \label{p:2}
If $G$ is a graph, then:
\begin{itemize}
\item $H^{(k)}(G,x) \ge 0$ for every $k\ge 0$, $x \in [0,\infty)$ and $d_u + d_v - 1 \geq k$,
\item $H(G,x) > 0$ on $(0,\infty)$ and  $H(G,x)$ is strictly increasing on $[0,\infty)$,
\item $H(G,x)$ is strictly convex on $[0,\infty)$ if and only if $G$ is not isomorphic to a union of path graphs $P_2$.
\end{itemize}
\end{proposition}

\begin{proof}
Since every coefficient of the polynomial $H(G,x)$ is non-negative, the first statement holds.

Since $\Deg_{\min} H(G,x)\ge 2\d-1 \ge 1$, we have $H(G,x) > 0$ and $H'(G,x)>0$ on $(0,\infty)$.

A graph $G$ is not isomorphic to a union of path graphs $P_2$ if and only if $d_u+d_v \ge 3$ for some edge $uv \in E(G)$;
this happens if and only if $G$ satisfies $\Deg H(G,x) \ge 2$; and this is equivalent to $H''(G,x)>0$ on $(0,\infty)$.
\end{proof}

Let us denote by $\mathfrak{G}$ the set of all regular and biregular connected graphs.
We say that a set of graphs $\{G_i\}_{i=1}^k$,
such that $G_i$ has maximum degree $\D_i$ and minimum degree $\d_i$ for each $1 \le i \le k$,
is \emph{coherent} if $G_i \subset \mathfrak{G}$ for every $1 \le i \le k$, and
$\D_i+\d_i = \D_j+\d_j$ for every $1 \le i,j \le k$.
We say that a graph is \emph{coherent} if the set of its connected components is coherent.

\smallskip

Given a graph $G$ and a vertex $v\in V(G)$, we denote by $N(v)$ the set of neighbors of $v$.

\begin{theorem} \label{t:p2}
Let $G$ be a graph. $x=0$ is the unique zero of $H(G,x)$ if and only if $G$ is coherent.
\end{theorem}

\begin{proof}
If $G$ is coherent, let us consider the set of its connected components $\{G_i\}_{i=1}^k$.
For each $1 \le i \le k$, $G_i$
is either a regular or a biregular graph with $m_i$ edges, maximum degree $\D_i$ and minimum degree $\d_i$;
hence, $H(G_i,x) = m_ix^{\D_i+\d_i-1}$.
So, $m=m_1+\dots +m_k$ is the cardinality of edges of $G$,
$H(G,x) = mx^{\D_1+\d_1-1}$
and $x=0$ is the unique zero of $H(G,x)$.

Assume now that $x=0$ is the unique zero of $H(G,x)$;
thus, $H(G,x) = ax^{b-1}$ for some positive integers $a,b,$ and $d_u+d_v=b$ for every $uv\in E(G)$.
Let us consider the set of connected components $\{G_i\}_{i=1}^k$ of $G$.
Fix $1 \le i \le k$, and denote by $\D_i$ and $\d_i$ the maximum and minimum degrees of $G$, respectively.
Thus, for each fixed vertex $u \in V(G_i)$ we have $d_v=b-d_u$ for every $uv\in E(G_i)$, and every $v \in N(u)$ has the same degree $b-d_u$.
In a similar way, if $w \in N(v)$, then $d_w=b-d_v=d_u$.
Since $G_i$ is a connected graph, $G_i$ is either regular (if $\D_i=\d_i$) or biregular (if $\D_i \neq \d_i$), and $G_i \subset \mathfrak{G}$.
Since $\D_i+\d_i = b$ for every $1 \le i \le k$, we conclude that $G$ is coherent.
\end{proof}

The following consequence of Theorem \ref{t:p2} shows that it is possible to characterize regular and biregular connected graphs in terms of the zeros of their harmonic polynomials.

\begin{corollary} \label{c:p2}
Let $G$ be a connected graph. $x=0$ is the unique zero of $H(G,x)$ if and only if $G$ is either a regular or a biregular graph.
\end{corollary}

The next result provides bounds of the harmonic index in terms of the values of the harmonic polynomial at the points $1$ and $1/2$.

\begin{proposition} \label{p:1bis}
If $G$ is a graph, then
$$
  H(G) \geq  2H(G,1/2),
%  \le H(G,1),
$$
and the equality in each inequality is attained if and only if $G$ is isomorphic to a union of path graphs $P_2$.
\end{proposition}

\begin{proof}
Hermite-Hadamard's inequality states that if $f:[0,1] \to \RR$ is a convex function, then
\begin{equation}
\label{eq:1bis}
f(1/2)
\le \int_0^1 f(x)\,dx,
%\le \frac{f(0)+f(1)}2 \,,
\end{equation}
and if $f$ is strictly convex, then the inequality is strict.

If $G$ is not isomorphic to a union of path graphs $P_2$, then Proposition \ref{p:2} gives
that $H(G,x)$ is a strictly convex function.
Thus, \eqref{eq:1bis} gives the result.
%$$
%2H(G,1/2) < H(G) = 2\int_0^1 H(G,x)\,dx < H(G,0)+H(G,1) = H(G,1).
%$$
%\indent
If $G$ is isomorphic to a union of $m$ path graphs $P_2$, then $H(G)=m$, $H(G,x)=mx$, $H(G,1/2)=m/2$.
%and $H(G,1)=m$.
Thus, $2H(G,1/2) = H(G)$.
\end{proof}

We say that a vertex $v\in V(G)$ in the graph $G$ is \emph{dominant} if $N(v) = V(G) \setminus \{v\}$.

\begin{proposition} \label{p:3}
Let $G$ be a graph with $n$ vertices, maximum degree $\D$ and minimum degree $\d$. Then:
\begin{itemize}
\item $x=0$ is a zero of $H(G,x)$ with multiplicity $\Deg_{\min} H(G,x)$, where $2\d-1 \le \Deg_{\min} H(G,x) \le \Deg H(G,x) \le 2\D-1$,
\item $\Deg H(G,x) \le 2n-3$, and $\Deg H(G,x) = 2n-3$ if and only if there are at least two dominant neighbors in $G$,
\item if $\G$ is a subgraph of $G$, then $\Deg H(\G,x) \le \Deg H(G,x)$ and \\
$\Deg_{\min} H(\G,x) \le \Deg_{\min} H(G,x)$.
\end{itemize}
\end{proposition}

\begin{proof}
Since
$$
H(G,x)=\sum_{j=\Deg_{\min} H(G,x)}^{\Deg H(G,x)} c_j \, x^j,
$$
for some constants $c_j$, $x=0$ is a zero of $H(G,x)$ with multiplicity $\Deg_{\min} H(G,x)$.
Since each $j$ in the previous sum can be written as $d_u+d_v-1$ for some $uv\in E(G)$, we have
$2\d-1 \le \Deg_{\min} H(G,x) \le \Deg H(G,x) \le 2\D-1$.

Since $\D \le n-1$, we have $\Deg H(G,x) \le 2n-3$.
We have $\Deg H(G,x) = 2n-3$ if and only if there is an edge $uv\in E(G)$ with $d_u=d_v= n-1$, and this holds is and only if
$u,v$ are dominant vertices in $G$.

Let $\G$ be a subgraph of $G$.
The last statement holds, since the degree of a vertex in $\G$ is at most its degree in $G$.
\end{proof}

The next result allows to obtain information about the connectedness, diameter and girth of a graph (the minimum length of its cycles)
in terms of the degree of its harmonic polynomial.

\begin{theorem} \label{t:p3}
Let $G$ be a graph with $n$ vertices.
If $\Deg H(G,x) \ge n$, then $g(G)=3$.
Furthermore, if $G$ is a triangle-free graph and $\Deg H(G,x) = n-1$, then $G$ is a connected graph and $\diam G \le 3$.
\end{theorem}

\begin{proof}
Since $g(G)=3$ if and only if $G$ is not triangle-free, it suffices to prove that if $G$ is a triangle-free graph, then $\Deg H(G,x) \le n-1$.
Since $G$ is a triangle-free graph, then $N(u)\cap N(v)=\emptyset$ for every $uv \in E(G)$.
Hence, $d_u+d_v\le n$ for every $uv \in E(G)$, and $\Deg H(G,x) \le n-1$.

Assume that $G$ is a triangle-free graph and $\Deg H(G,x) = n-1$.
Thus, there is an edge $uv\in E(G)$ with $d_u+d_v= n$.
Since $N(u)\cap N(v)=\emptyset$, we have $N(u)\cup N(v)= V(G)$ and $d(w,\{u,v\})\le 1$ for every $w \in V(G)$.
Consequently, $\diam G \le 3$ and $G$ is a connected graph.
\end{proof}

Denote by $K( p(x))$ the number of non-zero coefficients of the polynomial $p(x)$.

\begin{theorem} \label{t:p4}
Let $G$ be a graph with $m$ edges. Then:
\begin{itemize}
\item $1 \le K(H(G,x)) \le m$,
\item $K(H(G,x)) = 1$ if and only if $G$ is coherent,
\item $K(H(G,x)) = m$ if and only if $G$ is isomorphic to $P_2$.
\end{itemize}
\end{theorem}

\begin{proof}
%Since the sum in the definition of $H(G,x)$ has $m$ terms (with different exponent or not),
The first item is easy to see.

The proof of Theorem \ref{t:p2} gives that $G$ is coherent if and only if
$H(G,x) = ax^{b-1}$ for some positive integers $a,b$,
and this is equivalent to $K(H(G,x)) = 1$.

If $G$ is isomorphic to the path graph $P_2$, then it is a regular graph with just an edge, and the previous item gives $K(H(G,x)) = 1=m$.

Assume now that %$K(H(G,x)) = m$ and that $G$ is not isomorphic to $P_2$.
$G$ is not isomorphic to $P_2$.
We consider several cases.

$(1)$ $G$ is connected.
Thus, $3\le d_u+d_v \le m+1$ for every $uv \in E(G)$, i.e., $2\le d_u+d_v-1 \le m$.
Since the $m$ values of $d_u+d_v-1$ belong to a set of $m-1$ integers, there are two edges with the same value and
we conclude that $K(H(G,x)) \le m-1$.

$(2)$ $G$ is not connected.
So, $G$ has connected components $G_1,\dots,G_k,$ with $k \ge 2$.
Denote by $m_i$ the cardinality of the edges of $G_i$, thus $m = m_1+\dots+m_k$.

$(2.1)$ Assume that there exists some $1 \le j \le k$ such that $G_i$ is not isomorphic to $P_2$.
So, $(1)$ gives that $K(H(G_j,x)) \le m_j-1$, and this inequality and the first item give
$$
K(H(G,x))
\le \sum_{i=1}^{k} K(H(G_i,x))
\le \sum_{i=1}^{k} m_i - 1
= m-1.
$$
\indent
$(2.2)$ Assume that $G_i$ is isomorphic to $P_2$ for every $1 \le i \le k$.
So, $m =k \ge 2$,
$$
H(G,x)
= \sum_{i=1}^{m} H(G_i,x)
= \sum_{i=1}^{m} x
= mx,
$$
and $K(H(G,x)) = 1 \le k- 1 < m$.
\end{proof}

Theorem \ref{t:p4} has the following consequence.

\begin{corollary} \label{c:p4}
If $G$ is a graph with $m\ge 2$ edges, then $1 \le K(H(G,x)) \le m-1$.
\end{corollary}

\begin{proposition} \label{p:4}
Let $G$ be a graph with $n$ vertices, $m$ edges, maximum degree $\D$ and minimum degree $\d$. Then:
\begin{itemize}
\item $K(H(G,x)) \le \Deg H(G,x) - \Deg_{\min} H(G,x)+1$,
\item $K(H(G,x)) \le \min \{\,2\D - 2\d+1, \,m - 2\d+2 \}$,
\item if $G$ is a triangle-free graph, then $K(H(G,x)) \le n - 2\d+1$.
\end{itemize}
\end{proposition}

\begin{proof}
The first item holds since there are constants $c_j$ with
$$
H(G,x)=\sum_{j=\Deg_{\min} H(G,x)}^{\Deg H(G,x)} c_j \, x^j,
$$
\indent
The first item and the bounds in Proposition \ref{p:3} give
$$
K(H(G,x)) \le 2\D - 2\d+1 .
$$
Since $d_u+d_v \le m+1$ for every $uv \in E(G)$, we have $\Deg H(G,x) \le m$.
This inequality, the first item and the first item in Proposition \ref{p:3} give
$K(H(G,x)) \le m - 2\d+2$.

The third item is a consequence of the first one, the first item in Proposition \ref{p:3} and Theorem \ref{t:p3}.
\end{proof}

Given a graph $G$, we say that $\{d_u\}_{u \in V(G)}$
is the \emph{degree sequence} of $G$ (if $d_{v_1}=d_{v_2}$ for some $v_1,v_2 \in V(G)$, then the value $d_{v_1}=d_{v_2}$ appears just once in $\{d_u\}_{u \in V(G)}$).

\smallskip

Let us denote by $\lceil t \rceil$ the upper integer part of $t\in \RR$, i.e.,
the smallest integer greater or equal than $t$.

\begin{theorem} \label{t:sucdeg}
Let $G$ be a graph.
The following statements hold:
\begin{itemize}
\item if the degree sequence of $G$ has at most $r$ terms, then
$$
K(H(G,x)) \le \frac{r(r+1)}{2}\,,
$$
\item if $K(H(G,x)) \ge s$, then the degree sequence of $G$ has at least
$$
\left\lceil\frac{\sqrt{8s+1} -1}2 \right\rceil
$$
terms.
\end{itemize}
\end{theorem}

\begin{proof}
If the degree sequence of $G$ has at most $r$ terms, then the set of different values $d_u+d_v$
has cardinality at most $r(r+1)/2$
($2$-combinations with repetition of a set of $r$ elements).
Thus, $K(H(G,x)) \le r(r+1)/2$.

Assume that $K(H(G,x))=S\ge s$, and denote by $r$ the cardinality of the degree sequence of $G$.
The first item gives
$$
s \le S \le \frac{r(r+1)}{2}
\,,
\qquad
r^2+r-2s \ge 0,
\qquad
r \ge \frac{\sqrt{8s+1} -1}2
\,,
$$
and we obtain the desired inequality since $r$ is an integer.
\end{proof}

One can think that it might be possible to obtain a lower bound for $K(H(G,x))$
which is an increasing function of the cardinality of the degree sequence of $G$.
However, this is not possible, as the following result shows.

\begin{theorem} \label{t:sucdeg2}
Let $G$ be a connected graph with a degree sequence of cardinality $r$.
\begin{itemize}
\item If $r\le 2$, then $K(H(G,x))\ge 1$.
\item If $r>2$, then $K(H(G,x))\ge 2$.
\end{itemize}
Furthermore, the bounds are sharp for each $r$.
\end{theorem}

\begin{proof}
The first statement is a consequence of Theorem \ref{t:p4}.

Assume that $r > 2$.
Since $G$ is connected, there exist a path $\g=\{u_1,u_2,\dots,u_k\}$ in $G$ and three vertices in $V(G)\cap \g$ with different degrees.
Without loss of generality one can assume that
$$
d_{u_1} \notin  \{d_{u_2},\dots,d_{u_k}\} \quad \text{ and } \quad d_{u_k} \notin  \{d_{u_1},\dots,d_{u_{k-1}}\},
$$
since otherwise $u_1$ and/or $u_k$ can be removed from $\g$, and a shorter path with the same property is obtained.
Also, we can assume that $d_{u_2}=d_{u_{3}}=\dots = d_{u_{k-2}} = d_{u_{k-1}}$.
Thus, $d_{u_1} + d_{u_{2}}\neq d_{u_{2}} + d_{u_{k}} = d_{u_{k-1}} + d_{u_{k}}$ and,
since $u_1 u_{2}, u_{k-1} u_{k} \in E(G)$, we conclude $K(H(G,x))\ge 2$.

If $G$ is a star graph with $n$ vertices, then the degree sequence is $\{1,n-1\}$;
thus $r=1$ if $n=2$, and $r=2$ if $n>2$.
Since $H(G,x)= (n-1)\, x^{n-1}$, we have $K(H(G,x)) = 1$.

Consider the sequence $\{1,2,\dots,r\}$ with $r > 2$.
We are going to define a graph $T_r$ (in fact, $T_r$ is a tree) with degree sequence $\{1,2,\dots,r\}$ and $K(H(T_r,x)) = 2$.
Let us consider the (ordered) sequence $\{a_1,a_2,\dots,a_r\}$ obtained as a permutation of $\{1,2,\dots,r\}$ in the following way.
If $r$ is even, then
$$
\big\{a_1,a_2,\dots,a_r\big\}
=\Big\{\frac{r}2+1,\,\frac{r}2\,,\,\frac{r}2+2,\,\frac{r}2-1,\,\dots,r-1,\,2,\,r,\,1\Big\}.
$$
If $r$ is odd, then
$$
\big\{a_1,a_2,\dots,a_r\big\} =$$
$$=\Big\{\frac{r+1}2\,,\,\frac{r+1}2+1,\,\frac{r+1}2-1,\,\frac{r+1}2+2,\,\frac{r+1}2-2,\,\dots,r-1,\,2,\,r,\,1\Big\}.
$$
In both cases we have that $a_j+a_{j+1}$ is either $r+1$ or $r+2$ for each $1\le j < r$.
Consider a point $v_1$, which will be the root of $T_r$.
We define $T_r$ inductively on the distance $j$ from $v_1$.
We join $v_1$ with $a_1$ vertices (at distance $1$ from $v_1$).
If $u \in V(T_r)$ with $d_{T_r}(u,v_1)=j-1$ for some $1< j < r$, then we join $u$ with $a_j-1$ vertices (at distance $j$ from $v_1$).
Note that if $u \in V(T_r)$, then $d_{T_r}(u,v_1)=j-1$ for some $1\le j < r$ and $d_u = a_j$.
If $uv \in E(T_r)$, then without loss of generality we can assume that there exists $1\le j < r$ with $d_{T_r}(u,v_1)=j-1$ and $d_{T_r}(v,v_1)=j$.
Therefore, $d_u+d_v = a_j+a_{j+1}$ is either $r+1$ or $r+2$,
and so $K(H(T_r,x)) = 2$.
\end{proof}

\begin{theorem} \label{t:sucdeg3}
Let $G$ be a graph.
\begin{itemize}
\item If some connected component of $G$ has a degree sequence of cardinality $r > 2$, then $K(H(G,x))\ge 2$.
\item For each $r \ge 1$, there exists a graph with a degree sequence of cardinality $r$ and $K(H(G,x)) =1$.
\end{itemize}
\end{theorem}

\begin{proof}
If there is a connected component $G_i$ of $G$ with degree sequence of cardinality $r > 2$, then
Theorem \ref{t:sucdeg2} gives $K(H(G_i,x))\ge 2$, and  $K(H(G,x))\ge K(H(G_i,x))\ge 2$.

Fix any $r \ge 1$.

If $r$ is even, then define $G_r$ as the union of the complete bipartite graphs
$$
K_{1,r}, K_{2,r-1}, \dots, K_{r/2-1,r/2+2}, K_{r/2,r/2+1}.
$$

If $r$ is odd, then define $G_r$ as the union of the complete bipartite graphs
$$
K_{1,r}, K_{2,r-1}, \dots, K_{(r+1)/2-1,(r+1)/2+1}, K_{(r+1)/2,(r+1)/2}.
$$

In both cases, the degree sequence of $G_r$ has cardinality $r$.
If $m$ denotes the cardinality of $E(G_r)$, then $H(G_r,x) = m\,x^{r}$ and $K(H(G_r,x)) = 1$.
\end{proof}

Given a graph $G$, we say that the degree sequence of $G$ is even (respectively, odd) if $\{d_u\}_{u \in V(G)}$ is a subset of the even (respectively, odd) integers.

\begin{proposition} \label{p:sucdeg}
Let $G$ be a graph.
Then $H(G,x)$ is an odd function if and only if the degree sequence of each connected component of $G$ is either even or odd.
\end{proposition}

\begin{proof}
If the degree sequence of each connected component of $G$ is either even or odd, then $d_u+d_v-1$ is odd for every $uv \in E(G)$.
Since every exponent in $H(G,x)$ is odd, $H(G,x)$ is an odd function.

Assume now that $H(G,x)$ is an odd function.
Thus, $d_u+d_v$ is even for every $uv \in E(G)$.
Let us consider any fixed connected component $G_i$ of $G$.
If there is a vertex $u \in V(G_i)$ such that $d_u$ is even, then $d_v$ is even for every $v\in N(u)$.
Since $G_i$ is a connected graph, we conclude that the degree sequence of $G_i$ is even.
The same argument gives that if there is a vertex $u \in V(G_i)$ with $d_u$ odd, the degree sequence of $G_i$ is odd.
\end{proof}

We say that the graph $G$ has \emph{alternated degree} if $d_u$ and $d_v$ have different oddity for every $u,v \in V(G)$ with $uv \in E(G)$.

From the above definition, the following result is obtained.

\begin{proposition} \label{p:sucdeg2}
Let $G$ be a graph.
Then $H(G,x)$ is an even function if and only if $G$ has alternated degree.
\end{proposition}

%\begin{proof}
%If $G$ has alternated degree, then $d_u+d_v$ is odd for every $uv \in E(G)$.
%Thus, every exponent in $H(G,x)$ is even and $H(G,x)$ is an even function.
%
%Assume now that $H(G,x)$ is an even function.
%Hence, $d_u+d_v$ is odd for every $uv \in E(G)$, and $G$ has alternated degree.
%\end{proof}

An edge in a graph is said to be \emph{pendant} if one of its vertices has degree $1$.
A path with length two in a graph is said to be a \emph{pendant path} if it contains a pendant edge and a non-pendant edge.

\begin{proposition} \label{p:5}
Let $G$ be a graph.
Then, the cardinality of the pendant paths in $G$ is the coefficient of $x^2$ in $H(G,x)$.
\end{proposition}

\begin{proof}
There is a bijective correspondence between the pendant paths in $G$ and the edges $uv \in E(G)$ with $d_u=1$ and $d_v=2$ (i.e., $d_u+d_v-1=2$).
This gives the result.
\end{proof}

There are inequalities involving the harmonic and the first Zagreb indices (\cite{Ilic}, \cite[Theorem 2.5]{Xu}, \cite[p.234]{GFDMM}):

\begin{theorem} \label{t1}
Let $G$ be a graph with $m$ edges, maximum degree $\D$ and minimum degree $\d$. Then
$$
\dfrac{2m^{2}}{M_{1}(G)}
\le H(G)
\leq \dfrac{(\Delta+\delta)^{2}m^{2}}{2\Delta\delta M_{1}(G)} \,.
$$
The equality in the lower bound is attained if and only if $d_u+d_v$ is a constant for every $uv \in E(G)$.
The equality in the upper bound is attained if $G$ is regular.
\end{theorem}

We will use Theorem \ref{t1} in the proof of Proposition \ref{p:6} below.

\smallskip

Considering the Zagreb indices, Fath-Tabar \cite{FT} defined the first Zagreb polynomial as
$$
M_1(G, x) := \sum_{uv\in E(G)} x^{d_u + d_v}.
$$

The harmonic and the first Zagreb indices are related by Theorem \ref{t1}.
Moreover, the harmonic and the first Zagreb polynomials are related by the equality $M_1(G, x) = x\, H(G, x)$.

\smallskip

The next result provides more bounds of $\Deg_{\min} H(G,x)$ and $\Deg H(G,x)$.

\begin{proposition} \label{p:6}
Let $G$ be a graph with $n$ vertices, $m$ edges, maximum degree $\D$ and minimum degree $\d$. Then,
$$
2\d-1 \le \Deg_{\min} H(G,x) \le \frac{H'(G,1)}{m} \,,
\qquad
\frac{4m}{n}-1 \le \Deg H(G,x) \le 2\D-1.
%\frac{4m}{n}-1 \le \Deg_{\min} H(G,x) \le 2\D-1 \le 2n-3,
%\qquad
%2\d-1 \le \Deg H(G,x) \le \frac{H'(G,1)}{m} \,.
$$
\end{proposition}

\begin{proof}
The inequality $H(G) \leq n/2$ is a well-known upper bound for the harmonic index.
Theorem \ref{t1} gives the lower bound $H(G) \geq 2m^2/M_1(G)$.
Given $j\in \NN$, let us define $c_j=c_j(G)$ as the cardinality of the set $\{uv\in E(G)\, |\; d_u+d_v-1=j\}$.
We can write
$$
H(G, x)
= \sum_{j=\Deg_{\min} H(G,x)}^{\Deg H(G,x)} c_j x^{j} ,
\qquad
\text{with }
\;
\sum_{j=\Deg_{\min} H(G,x)}^{\Deg H(G,x)} c_j =m.
$$
Thus, we have
$$
\begin{aligned}
\frac{n}2
& \ge H(G)
= 2\int^1_0 H(G, x) \,dx
= \sum_{j=\Deg_{\min} H(G,x)}^{\Deg H(G,x)} \frac{2c_j}{j+1}
\\
& \ge \sum_{j=\Deg_{\min} H(G,x)}^{\Deg H(G,x)} \frac{2c_j}{\Deg H(G,x)+1}
= \frac{2m}{\Deg H(G,x)+1}\,,
\\
\Deg H(G,x)
& \ge \frac{4m}{n}-1,
\\
%\frac{2m}{m+1}
%& \le H(G)
%= \sum_{j=\Deg_{\min} H(G,x)}^{\Deg H(G,x)} \frac{2c_j}{j+1}
%\le \frac{2m}{\Deg_{\min} H(G,x)+1}\,,
%\\
\frac{2m^2}{M_1(G)}
& \le H(G)
= \sum_{j=\Deg_{\min} H(G,x)}^{\Deg H(G,x)} \frac{2c_j}{j+1}
\le \frac{2m}{\Deg_{\min} H(G,x)+1}\,,
\\
\Deg_{\min} H(G,x)
& \le \frac{M_1(G)}{m}-1
= \frac{M_1(G)-m}{m}
= \frac{H'(G,1)}{m}
\,.
\end{aligned}
$$
Proposition \ref{p:3} provides the other inequalities.
\end{proof}

The next result allows to bound the harmonic index of a graph by using several parameters of its harmonic polynomial.

\smallskip

Given a graph $G$, let us denote by $c_{\min}(G)$ and $c_{\max}(G)$ the coefficients of
$x^{\Deg_{\min} H(G,x)}$ and $x^{\Deg H(G,x)}$ in $H(G,x)$,
respectively.

\begin{proposition} \label{p:7}
Let $G$ be a graph with $m$ edges. Then,
$$
\frac{2c_{\min}(G)}{\Deg_{\min} H(G,x)+1} + \frac{2m-2c_{\min}(G)}{\Deg H(G,x)+1}
\le H(G) 
\le $$
$$\le \frac{2c_{\max}(G)}{\Deg H(G,x)+1} + \frac{2m-2c_{\max}(G)}{\Deg_{\min} H(G,x)+1}
\,.
$$
\end{proposition}

\begin{proof}
As in the proof of Proposition \ref{p:6}, we obtain
$$
H(G)
= 2\int^1_0 H(G, x) \,dx
= \sum_{j=\Deg_{\min} H(G,x)}^{\Deg H(G,x)} \frac{2c_j}{j+1}
\,.
$$
Hence,
$$
\begin{aligned}
H(G)
& = \frac{2c_{\min}(G)}{\Deg_{\min} H(G,x)+1} + \sum_{j=\Deg_{\min} H(G,x)+1}^{\Deg H(G,x)} \frac{2c_j}{j+1}
\\
& \ge \frac{2c_{\min}(G)}{\Deg_{\min} H(G,x)+1} + \sum_{j=\Deg_{\min} H(G,x)+1}^{\Deg H(G,x)} \frac{2c_j}{\Deg H(G,x)+1}
\\
& = \frac{2c_{\min}(G)}{\Deg_{\min} H(G,x)+1} + \frac{2m-2c_{\min}(G)}{\Deg H(G,x)+1}
\,,
\\
H(G)
& = \frac{2c_{\max}(G)}{\Deg H(G,x)+1} + \sum_{j=\Deg_{\min} H(G,x)}^{\Deg H(G,x)-1} \frac{2c_j}{j+1}
\\
& \le \frac{2c_{\max}(G)}{\Deg H(G,x)+1} + \sum_{j=\Deg_{\min} H(G,x)}^{\Deg H(G,x)-1} \frac{2c_j}{\Deg_{\min} H(G,x)+1}
\\
& = \frac{2c_{\max}(G)}{\Deg H(G,x)+1} + \frac{2m-2c_{\max}(G)}{\Deg_{\min} H(G,x)+1}
\,.
\end{aligned}
$$
\end{proof}

Although two non-isomorphic graphs can have the same harmonic polynomial,
Theorem \ref{t:mu} below shows that two graphs with the same harmonic polynomial
have to be similar.

\smallskip

For each function $\mu: \NN \rightarrow (0,\infty)$, let us define its associated topological indices
$$
T_\mu(G) = \sum_{uv\in E(G)}\mu(d_u + d_v)\,,
\qquad
U_\mu(G) = \prod_{uv\in E(G)}\mu(d_u + d_v)\,.
$$
In particular, if $\mu(t)=t^\a$, then $T_\mu=\chi_{\a}$.
The \emph{modified first multiplicative Zagreb index} is defined in \cite{EIG} by
$\Pi_1^*(G) = \prod_{uv\in E(G)}(d_u + d_v)\,.$
In particular, if $\mu(t)=t$, then $U_\mu=\Pi_1^*$.

\begin{theorem} \label{t:mu}
If two graphs $G_1$ and $G_2$ have the same harmonic polynomial, then
$T_\mu(G_1) =T_\mu(G_2)$ and $U_\mu(G_1) =U_\mu(G_2)$ for every function $\mu: \NN \rightarrow (0,\infty)$.
In particular, $\chi_{\a}(G_1) =\chi_{\a}(G_2)$ for every $\a \in \RR$, and
$\Pi_1^*(G_1) =\Pi_1^*(G_2)$.
\end{theorem}

\begin{proof}
As in the proof of Proposition \ref{p:6},
given a graph $G$ and $j\in \NN$, we define $c_j(G)$ as the cardinality of the set $\{uv\in E(G)\, |\; d_u+d_v-1=j\}$.
Thus, $H(G,x)=\sum_j c_j(G)\,x^{j}$.
If $H(G_1,x)=H(G_2,x)$, then $c_j(G_1) = c_j(G_2)$ for every $j\in \NN$.
Since $T_\mu(G)=\sum_j c_j(G)\mu(j+1)$ and $U_\mu(G) = \prod_{j}\mu(j+1)^{c_j(G)}$
for every function $\mu: \NN \rightarrow (0,\infty)$, we conclude that $T_\mu(G_1) =T_\mu(G_2)$ and $U_\mu(G_1) =U_\mu(G_2)$.
\end{proof}

%\begin{remark} \label{r:mu}
We want to remark that if we consider a function $\mu: \NN \rightarrow \CC$ in the definition $T_\mu$,
then the argument in the proof of Theorem \ref{t:mu} also works.
Thus, we can consider a family of functions $\{\mu_z\}$, where $z$ is a complex variable,
and we can define for each graph $G$ the complex function $F_{G}(z):=T_{\mu_z}(G)$.
So, if two graphs $G_1$ and $G_2$ have the same harmonic polynomial, then the complex functions
$F_{G_1}(z)$ and $F_{G_1}(z)$ are the same.
This holds, in particular, for the holomorphic function
$F_{G}(z):= \sum_{uv\in E(G)}(d_u + d_v)^z$.
%\end{remark}

\smallskip

\section*{Acknowledgements}
%Supported in part by two grants from Ministerio de Economía y Competitividad, Agencia Estatal de Investigación (AEI) and Fondo Europeo de Desarrollo Regional (FEDER) (MTM2016-78227-C2-1-P and MTM2017-90584-REDT), Spain.
We would like to thank the referees for their careful reading of the manuscript and for some helpful suggestions that have improved the paper.

\end{document}